\documentstyle[12pt,epsf]{amlt}
\begin{document}


\startingpage{1}



\def\descriptionlabel#1{\normalshape\rm\kern-.5em#1}

\newcommand{\mycirc}{\hbox{\raise.25ex\hbox{${\scriptstyle
\circ}$}}\ }
\newcommand{\ig}{\ell}
\newcommand{\ttot}{\tau_{\mbox{\scriptsize total}}}
\renewcommand{\varphi}{\wp}
\newcommand{\R}{\bold{R}} 
\newcommand{\h}{\bold{H}} 
\newcommand{\s}{\bold{S}} 
\newcommand{\E}{\bold{E}} 
\newcommand{\p}{{\cal P}}
\newcommand{\m}{{\cal M}}
\newcommand{\thefig}{\thefigure}
\newcommand{\bd}[1]{{\bf #1}}  
\newcommand{\fn}[1]{{\tiny #1}}
\newcommand{\cp}[1]{{\sc #1}}  
\newcommand{\ic}[1]{{\it #1}}  
\newcommand{\es}[1]{{\em #1}}  
\newcommand{\ro}[1]{\mathrm{#1}}  

\newcommand{\Kext}{K_{\mbox{\scriptsize ext}}}
\def\QED{\enddemo}
\newcommand{\sphi}{$\bigcirc_\infty$}  
\newcommand{\ij}{\ro{InjRad}}  
\let\mcol=\multicolumn
\def\0{\kern.5em}		
\newcommand{\po}{{\cal P}}  
\newcommand{\poi}{{\cal P}^{- 1}}  
\def\di{\mathop {\rm dist}\nolimits}     
\def\det{\mathop {\rm det}\nolimits}     
\def\tr{\mathop {\rm tr}\nolimits}     
\def\grad{\mathop {\rm grad}\nolimits}     
\newcommand{\len}{\ell}         
\def\arctanh{\mathop {\rm arctanh}\nolimits}
\def\arcsec{\mathop {\rm arcsec}\nolimits}
\def\sign{\mathop {\rm sign}\nolimits}
\def\Aut{\mathop {\rm Aut}\nolimits}
\def\Isom{\mathop {\mbox{\rm Isom}}\nolimits}
\newcommand{\degree}{^{\circ}}
\newcommand{\degrees}{^{\circ}}
\newcommand{\arccosh}{\ro{arccosh}}
\def\arccosh{\mathop {\rm arccosh}\nolimits}
\def\arcsinh{\mathop {\rm arcsinh}\nolimits}
\def\diam{\mathop {\rm diam}\nolimits}
\def\id{\mathop {\rm id}\nolimits}
\newcommand{\Bd}{{\partial }} 
\newcommand{\ip}[3]{\langle #1 \rangle_{#2 , #3}}   
\newcommand{\nr}[3]{\| #1 \|_{#2 , #3}}  
\newcommand{\FIG}[1]{\vspace{ #1 }}   

\newcommand{\leb}[1]{} 
\newcommand{\lab}[1]{\label{#1}}
\newcommand{\rf}[1]{\ref{#1}}
\newcommand{\numb}[1]{(\ref{#1})}
\newcommand{\thm}[2]{{\proclaim{#1} #2 \endproclaim}}
\newcommand{\ie}[1]{{\index{#1}}}
\newcommand{\mpf}[1]{\marginpar{\fn#1}}


\def\sun{{\hbox{$\odot$}}}
\def\la{\mathrel{\mathchoice {\vcenter{\offinterlineskip\halign{\hfil
$\displaystyle##$\hfil\cr<\cr\noalign{\vskip1.5pt}\sim\cr}}}
{\vcenter{\offinterlineskip\halign{\hfil$\textstyle##$\hfil\cr<\cr
\noalign{\vskip1.0pt}\sim\cr}}}
{\vcenter{\offinterlineskip\halign{\hfil$\scriptstyle##$\hfil\cr<\cr
\noalign{\vskip0.5pt}\sim\cr}}}
{\vcenter{\offinterlineskip\halign{\hfil$\scriptscriptstyle##$\hfil
\cr<\cr\noalign{\vskip0.5pt}\sim\cr}}}}}
\def\ga{\mathrel{\mathchoice {\vcenter{\offinterlineskip\halign{\hfil
$\displaystyle##$\hfil\cr>\cr\noalign{\vskip1.5pt}\sim\cr}}}
{\vcenter{\offinterlineskip\halign{\hfil$\textstyle##$\hfil\cr>\cr
\noalign{\vskip1.0pt}\sim\cr}}}
{\vcenter{\offinterlineskip\halign{\hfil$\scriptstyle##$\hfil\cr>\cr
\noalign{\vskip0.5pt}\sim\cr}}}
{\vcenter{\offinterlineskip\halign{\hfil$\scriptscriptstyle##$\hfil
\cr>\cr\noalign{\vskip0.5pt}\sim\cr}}}}}
\def\sq{\hbox{\rlap{$\sqcap$}$\sqcup$}}
\def\degr{\hbox{$^\circ$}}
\def\arcmin{\hbox{$^\prime$}}
\def\arcsec{\hbox{$^{\prime\prime}$}}
\def\utw{\smash{\rlap{\lower5pt\hbox{$\sim$}}}}
\def\udtw{\smash{\rlap{\lower6pt\hbox{$\approx$}}}}
\def\fd{\hbox{$.\!\!^{\rm d}$}}
\def\fh{\hbox{$.\!\!^{\rm h}$}}
\def\fm{\hbox{$.\!\!^{\rm m}$}}
\def\fs{\hbox{$.\!\!^{\rm s}$}}
\def\fdg{\hbox{$.\!\!^\circ$}}
\def\farcm{\hbox{$.\mkern-4mu^\prime$}}
\def\farcs{\hbox{$.\!\!^{\prime\prime}$}}
\def\fp{\hbox{$.\!\!^{\scriptscriptstyle\rm p}$}}
\def\getsto{\mathrel{\mathchoice {\vcenter{\offinterlineskip
\halign{\hfil$\displaystyle##$\hfil\cr\gets\cr\to\cr}}}
{\vcenter{\offinterlineskip\halign{\hfil$\textstyle##$\hfil\cr
\gets\cr\to\cr}}}
{\vcenter{\offinterlineskip\halign{\hfil$\scriptstyle##$\hfil\cr
\gets\cr\to\cr}}}
{\vcenter{\offinterlineskip\halign{\hfil$\scriptscriptstyle##$\hfil\cr
\gets\cr\to\cr}}}}}
\def\cor{\mathrel{\mathchoice {\hbox{$\widehat=$}}{\hbox{$\widehat=$}}
{\hbox{$\scriptstyle\hat=$}}
{\hbox{$\scriptscriptstyle\hat=$}}}}
\def\grole{\mathrel{\mathchoice {\vcenter{\offinterlineskip\halign{\hfil
$\displaystyle##$\hfil\cr>\cr\noalign{\vskip-1.5pt}<\cr}}}
{\vcenter{\offinterlineskip\halign{\hfil$\textstyle##$\hfil\cr
>\cr\noalign{\vskip-1.5pt}<\cr}}}
{\vcenter{\offinterlineskip\halign{\hfil$\scriptstyle##$\hfil\cr
>\cr\noalign{\vskip-1pt}<\cr}}}
{\vcenter{\offinterlineskip\halign{\hfil$\scriptscriptstyle##$\hfil\cr
>\cr\noalign{\vskip-0.5pt}<\cr}}}}}
\def\lid{\mathrel{\mathchoice {\vcenter{\offinterlineskip\halign{\hfil
$\displaystyle##$\hfil\cr<\cr\noalign{\vskip1.5pt}=\cr}}}
{\vcenter{\offinterlineskip\halign{\hfil$\textstyle##$\hfil\cr<\cr
\noalign{\vskip1pt}=\cr}}}
{\vcenter{\offinterlineskip\halign{\hfil$\scriptstyle##$\hfil\cr<\cr
\noalign{\vskip0.5pt}=\cr}}}
{\vcenter{\offinterlineskip\halign{\hfil$\scriptscriptstyle##$\hfil\cr
<\cr\noalign{\vskip0.5pt}=\cr}}}}}
\def\gid{\mathrel{\mathchoice {\vcenter{\offinterlineskip\halign{\hfil
$\displaystyle##$\hfil\cr>\cr\noalign{\vskip1.5pt}=\cr}}}
{\vcenter{\offinterlineskip\halign{\hfil$\textstyle##$\hfil\cr>\cr
\noalign{\vskip1pt}=\cr}}}
{\vcenter{\offinterlineskip\halign{\hfil$\scriptstyle##$\hfil\cr>\cr
\noalign{\vskip0.5pt}=\cr}}}
{\vcenter{\offinterlineskip\halign{\hfil$\scriptscriptstyle##$\hfil\cr
>\cr\noalign{\vskip0.5pt}=\cr}}}}}
\def\sol{\mathrel{\mathchoice {\vcenter{\offinterlineskip\halign{\hfil
$\displaystyle##$\hfil\cr\sim\cr\noalign{\vskip-0.2mm}<\cr}}}
{\vcenter{\offinterlineskip
\halign{\hfil$\textstyle##$\hfil\cr\sim\cr<\cr}}}
{\vcenter{\offinterlineskip
\halign{\hfil$\scriptstyle##$\hfil\cr\sim\cr<\cr}}}
{\vcenter{\offinterlineskip
\halign{\hfil$\scriptscriptstyle##$\hfil\cr\sim\cr<\cr}}}}}
\def\sog{\mathrel{\mathchoice {\vcenter{\offinterlineskip\halign{\hfil
$\displaystyle##$\hfil\cr\sim\cr\noalign{\vskip-0.2mm}>\cr}}}
{\vcenter{\offinterlineskip
\halign{\hfil$\textstyle##$\hfil\cr\sim\cr>\cr}}}
{\vcenter{\offinterlineskip
\halign{\hfil$\scriptstyle##$\hfil\cr\sim\cr>\cr}}}
{\vcenter{\offinterlineskip
\halign{\hfil$\scriptscriptstyle##$\hfil\cr\sim\cr>\cr}}}}}
\def\lse{\mathrel{\mathchoice {\vcenter{\offinterlineskip\halign{\hfil
$\displaystyle##$\hfil\cr<\cr\noalign{\vskip1.5pt}\simeq\cr}}}
{\vcenter{\offinterlineskip\halign{\hfil$\textstyle##$\hfil\cr<\cr
\noalign{\vskip1pt}\simeq\cr}}}
{\vcenter{\offinterlineskip\halign{\hfil$\scriptstyle##$\hfil\cr<\cr
\noalign{\vskip0.5pt}\simeq\cr}}}
{\vcenter{\offinterlineskip
\halign{\hfil$\scriptscriptstyle##$\hfil\cr<\cr
\noalign{\vskip0.5pt}\simeq\cr}}}}}
\def\gse{\mathrel{\mathchoice {\vcenter{\offinterlineskip\halign{\hfil
$\displaystyle##$\hfil\cr>\cr\noalign{\vskip1.5pt}\simeq\cr}}}
{\vcenter{\offinterlineskip\halign{\hfil$\textstyle##$\hfil\cr>\cr
\noalign{\vskip1.0pt}\simeq\cr}}}
{\vcenter{\offinterlineskip\halign{\hfil$\scriptstyle##$\hfil\cr>\cr
\noalign{\vskip0.5pt}\simeq\cr}}}
{\vcenter{\offinterlineskip
\halign{\hfil$\scriptscriptstyle##$\hfil\cr>\cr
\noalign{\vskip0.5pt}\simeq\cr}}}}}
\def\grole{\mathrel{\mathchoice {\vcenter{\offinterlineskip\halign{\hfil
$\displaystyle##$\hfil\cr>\cr\noalign{\vskip-1.5pt}<\cr}}}
{\vcenter{\offinterlineskip\halign{\hfil$\textstyle##$\hfil\cr
>\cr\noalign{\vskip-1.5pt}<\cr}}}
{\vcenter{\offinterlineskip\halign{\hfil$\scriptstyle##$\hfil\cr
>\cr\noalign{\vskip-1pt}<\cr}}}
{\vcenter{\offinterlineskip\halign{\hfil$\scriptscriptstyle##$\hfil\cr
>\cr\noalign{\vskip-0.5pt}<\cr}}}}}
\def\leogr{\mathrel{\mathchoice {\vcenter{\offinterlineskip\halign{\hfil
$\displaystyle##$\hfil\cr<\cr\noalign{\vskip-1.5pt}>\cr}}}
{\vcenter{\offinterlineskip\halign{\hfil$\textstyle##$\hfil\cr
<\cr\noalign{\vskip-1.5pt}>\cr}}}
{\vcenter{\offinterlineskip\halign{\hfil$\scriptstyle##$\hfil\cr
<\cr\noalign{\vskip-1pt}>\cr}}}
{\vcenter{\offinterlineskip\halign{\hfil$\scriptscriptstyle##$\hfil\cr
<\cr\noalign{\vskip-0.5pt}>\cr}}}}}
\def\loa{\mathrel{\mathchoice {\vcenter{\offinterlineskip\halign{\hfil
$\displaystyle##$\hfil\cr<\cr\noalign{\vskip1.5pt}\approx\cr}}}
{\vcenter{\offinterlineskip\halign{\hfil$\textstyle##$\hfil\cr<\cr
\noalign{\vskip1.0pt}\approx\cr}}}
{\vcenter{\offinterlineskip\halign{\hfil$\scriptstyle##$\hfil\cr<\cr
\noalign{\vskip0.5pt}\approx\cr}}}
{\vcenter{\offinterlineskip\halign{\hfil$\scriptscriptstyle##$\hfil\cr
<\cr\noalign{\vskip0.5pt}\approx\cr}}}}}
\def\goa{\mathrel{\mathchoice {\vcenter{\offinterlineskip\halign{\hfil
$\displaystyle##$\hfil\cr>\cr\noalign{\vskip1.5pt}\approx\cr}}}
{\vcenter{\offinterlineskip\halign{\hfil$\textstyle##$\hfil\cr>\cr
\noalign{\vskip1.0pt}\approx\cr}}}
{\vcenter{\offinterlineskip\halign{\hfil$\scriptstyle##$\hfil\cr>\cr
\noalign{\vskip0.5pt}\approx\cr}}}
{\vcenter{\offinterlineskip\halign{\hfil$\scriptscriptstyle##$\hfil\cr
>\cr\noalign{\vskip0.5pt}\approx\cr}}}}}
\def\bbbr{\mathbold{R}} 
\def\bbbm{\mathbold{M}}
\def\bbbh{\mathbold{H}}
\def\bbbd{\mathbold{D}}
\def\bbbone{{\mathchoice {\rm 1\mskip-4mu l} {\rm 1\mskip-4mu l}
{\rm 1\mskip-4.5mu l} {\rm 1\mskip-5mu l}}}
\def\bbbc{{\mathchoice {\setbox0=\hbox{$\displaystyle\rm C$}\hbox{\hbox
to0pt{\kern0.4\wd0\vrule height0.9\ht0\hss}\box0}}
{\setbox0=\hbox{$\textstyle\rm C$}\hbox{\hbox
to0pt{\kern0.4\wd0\vrule height0.9\ht0\hss}\box0}}
{\setbox0=\hbox{$\scriptstyle\rm C$}\hbox{\hbox
to0pt{\kern0.4\wd0\vrule height0.9\ht0\hss}\box0}}
{\setbox0=\hbox{$\scriptscriptstyle\rm C$}\hbox{\hbox
to0pt{\kern0.4\wd0\vrule height0.9\ht0\hss}\box0}}}}
\def\bbbe{{\bf E}}
\def\bbbq{{\mathchoice {\setbox0=\hbox{$\displaystyle\rm Q$}\hbox{\raise
0.15\ht0\hbox to0pt{\kern0.4\wd0\vrule height0.8\ht0\hss}\box0}}
{\setbox0=\hbox{$\textstyle\rm Q$}\hbox{\raise
0.15\ht0\hbox to0pt{\kern0.4\wd0\vrule height0.8\ht0\hss}\box0}}
{\setbox0=\hbox{$\scriptstyle\rm Q$}\hbox{\raise
0.15\ht0\hbox to0pt{\kern0.4\wd0\vrule height0.7\ht0\hss}\box0}}
{\setbox0=\hbox{$\scriptscriptstyle\rm Q$}\hbox{\raise
0.15\ht0\hbox to0pt{\kern0.4\wd0\vrule height0.7\ht0\hss}\box0}}}}
\def\bbbt{{\mathchoice {\setbox0=\hbox{$\displaystyle\rm
T$}\hbox{\hbox to0pt{\kern0.3\wd0\vrule height0.9\ht0\hss}\box0}}
{\setbox0=\hbox{$\textstyle\rm T$}\hbox{\hbox
to0pt{\kern0.3\wd0\vrule height0.9\ht0\hss}\box0}}
{\setbox0=\hbox{$\scriptstyle\rm T$}\hbox{\hbox
to0pt{\kern0.3\wd0\vrule height0.9\ht0\hss}\box0}}
{\setbox0=\hbox{$\scriptscriptstyle\rm T$}\hbox{\hbox
to0pt{\kern0.3\wd0\vrule height0.9\ht0\hss}\box0}}}}
\def\bbbs{{\mathchoice
{\setbox0=\hbox{$\displaystyle     \rm S$}\hbox{\raise0.5\ht0\hbox
to0pt{\kern0.35\wd0\vrule height0.45\ht0\hss}\hbox
to0pt{\kern0.55\wd0\vrule height0.5\ht0\hss}\box0}}
{\setbox0=\hbox{$\textstyle        \rm S$}\hbox{\raise0.5\ht0\hbox
to0pt{\kern0.35\wd0\vrule height0.45\ht0\hss}\hbox
to0pt{\kern0.55\wd0\vrule height0.5\ht0\hss}\box0}}
{\setbox0=\hbox{$\scriptstyle      \rm S$}\hbox{\raise0.5\ht0\hbox
to0pt{\kern0.35\wd0\vrule height0.45\ht0\hss}\raise0.05\ht0\hbox
to0pt{\kern0.5\wd0\vrule height0.45\ht0\hss}\box0}}
{\setbox0=\hbox{$\scriptscriptstyle\rm S$}\hbox{\raise0.5\ht0\hbox
to0pt{\kern0.4\wd0\vrule height0.45\ht0\hss}\raise0.05\ht0\hbox
to0pt{\kern0.55\wd0\vrule height0.45\ht0\hss}\box0}}}}

%
%

\def\bbbz{{\mathchoice {\hbox{$\sf\textstyle Z\kern-0.4em Z$}}
{\hbox{$\sf\textstyle Z\kern-0.4em Z$}}
{\hbox{$\sf\scriptstyle Z\kern-0.3em Z$}}
{\hbox{$\sf\scriptscriptstyle Z\kern-0.2em Z$}}}}

\def\diameter{{\ifmmode\oslash\else$\oslash$\fi}}

\newcommand{\sh}{{\euf s}}
\newcommand{\cbar}{\overline{\bbc}}
\newcommand{\crat}{{\euf c}}
\newcommand{\ccirc}[1]{{\cal C}_{#1}}
\newcommand{\hyp}{\bbbh}
\newcommand{\htwo}{\hyp^2}
\newcommand{\hthree}{\hyp^3}
\newcommand{\etwo}{\bbbe^2}
\newcommand{\ethree}{\bbbe^3}
\newcommand{\reals}{\bbbr}
\newcommand{\comps}{\bbbc}
\def\card{\mathop {\rm card}\nolimits}
\def\grad{\mathop {\rm grad}\nolimits}
\newtheorem{Theorem}{Theorem}[section]
\renewcommand{\theTheorem}{\arabic{section}.\arabic{Theorem}}
\newtheorem{Lemma}[Theorem]{Lemma}
\newtheorem{Remark}[Theorem]{Remark}
\newtheorem{Example}[Theorem]{Example}
\newtheorem{Observation}[Theorem]{Observation}
\newtheorem{Proposition}[Theorem]{Proposition}
\newtheorem{Corollary}[Theorem]{Corollary}
\newtheorem{Definition}[Theorem]{Definition}
\newtheorem{Notes}{Notes}
\renewcommand{\theNotes}{}
\newtheorem{Fact}{Fact}
\renewcommand{\theFact}{}
\newtheorem{ParTheorem}{Parametrization Theorem}
\renewcommand{\theParTheorem}{}
\newtheorem{Note}{Note}
\renewcommand{\theNote}{}


\title{A norm on homology of surfaces and counting simple geodesics}
\shorttitle{Geodesics and homology} 

\acknowledgements{
Igor Rivin has been supported in part by the
Institute for Advanced Study under the NSF Grant  \#DMS-9304580. Igor
Rivin would like to thank Caltech for its hospitality during
the preparation of this paper}
\twoauthors{Greg McShane}{Igor Rivin}
\institutions{Ecole Normale Superieure, Lyon\\
Melbourne University, Victoria}

\intro
Let $S$ be a hyperbolic surface of finite volume, and define $N_S(L)$
to be the number of simple (that is, without self-intersections)
closed geodesics on $S$ whose length does not exceed $L.$ 
This quantity arises naturally in a number of contexts, and hence it
is useful to obtain some estimates on the order of growth of $N_S(L)$
as $L$ grows large. Although a number of people have worked on this
question, the answers have not, on the whole, been completely
satisfactory. For general cusped $S$ Birman and Series \cite{birser}
have shown that the $N_S(L)$ grows at most polynomially as a function
of $L.$ They bound the degree of the polynomial by a function
of the genus, and note that the bound is very far from optimal. For
the simplest hyperbolic surface -- the once punctured torus,
the best result to date has been that of Beardon, Lehner, and
Sheingorn \cite{bls}, who showed that there $N_S(L)$ grows at least
linearly, and at most quadratically in $L.$ The case of the punctured
torus is interesting for reasons other than its relative simplicity --
simple geodesics on the punctured torus have deep connections with
diophantine approximation, see for example \cite{haas}. 

The purpose of this note is to prove the following estimate:

\thm{Theorem}
{\label{maintheorem}
Let $S$ be a punctured torus equipped with a complete hyperbolic
metric of finite volume. Then, $\lim_{L\rightarrow \infty}
N_S(L) = c_S (L^2) + O(L \log L),$ where $c_S$ depends on the
hyperbolic metric. 
}

The proof of Theorem \ref{maintheorem} is based on
some preliminaries from hyperbolic geometry (Section \ref{poincare}),
topology of surfaces (Section \ref{nielsen}), ),the general metric
theory of groups acting on spaces (section \ref{cannon}). We final
ingredient is a norm on the first homology of surfaces coming from the
study of of curve systems (Section \ref{multic}, and subsequent).

\demo{Note} An estimate of type \ref{maintheorem} for the modular
torus (but stated in the language of Markoff spectra) is obtained by
Don Zagier in \ref{zagier}. Zagier's methods are computational. Geoff
Mess told the authors that he has  independently discovered the norm
on homology used here. Mess' result has not been published. 
\enddemo

\section{Simple geodesics do not enter cusps}
\label{poincare}
The result of this section goes back to Poincar\'e, while the sharp
version for the torus stated below in Theorem \ref{mcshane} is proved
in Greg McShane's 1991 Warwick thesis \cite{mcthes}. First, a
definition:

\demo{Definition} A {\em cusp region} is a neighborhood of the cusp of
$S,$ bounded by a horocycle.
\enddemo

\thm{Theorem}
{\label{mcshane}
Let $\epsilon > 0.$ Any punctured torus has a cusp region with
bounding curve of length $4 - \epsilon,$ and this bound is optimal. No
simple closed geodesic intersects a cusp region with boundary curve of
length $4 - \epsilon.$
}

For current purposes it will be sufficient to prove the following very
easy theorem, weaker than Theorem \ref{mcshane} in the case of a punctured
torus.

\thm{Theorem}
{\label{mcsimple}
Any cusped hyperbolic surface $S$  has a cusp region with
bounding curve of length $2.$ No
simple closed geodesic intersects a cusp region with boundary curve of
length $2.$
}

\demo{Proof}  In the  upper half-plane model of $H^2,$ consider a
fundamental domain of $S,$ arranged in such a way that the parabolic
element preserving the cusp in question is $\lambda:z\rightarrow z+1.$
A simple closed geodesic $g$ of $S$ lifts to a geodesic in $H^2,$
which is represented by a semicircle $\tilde{g}$ in the half-space
model. Since $g$ is simple, it follows that $\tilde{g}\cap
\lambda(\tilde{g}) = \emptyset.$ Thus, the radius of the semicircle
representing $\tilde{g}$ is smaller than $\frac12.$ The same goes for any
non-vertical boundary component of the fundamental domain of $S,$ and
thus the horocyclic arc joining (for example) $\frac12 i$ and $1+\frac12 i$ is
entirely contained in the  fundamental domain of $S,$ and meets no
lift of any closed geodesic of $S.$ This horocyclic arc (which
projects to a closed horocycle in $S$ has length 2.
\enddemo

\section{Primitive elements in the fundamental group of the punctured
torus}
\label{nielsen}
The reader should recall that the fundamental group $\pi_1(T)$ of the
punctured torus is the free group on two generators $F_2 = \langle
s, t \rangle,$ where $s$ and $t$ are the standard generators. 

\demo{Definition} An element $\gamma\in F_2$ is called primitive, if
there exists an automorphism $\phi$ of $F_2,$ such that $\psi(\gamma)=
s.$ If $\phi(\delta)=t,$ then $\gamma$ and $\delta$ are called
{\em associated primitives.} 
\enddemo

\thm{Fact}
{\label{fact1}
The outer automorphism group of $F_2$ is isomorphic to the mapping
class group of the punctured torus.
}

This was apparently first observed by Max Dehn,
though Fact \ref{fact1} has now passed into the folklore. One argument
proceeds roughly as follows:

The outer automorphism group is generated by the so-called Nielsen
transformation, which either permute the basis $x, y,$ or transform it into 
$x y, y,$ or $x y^{-1}, y.$ In the case of a punctured torus these
transformations can be topologically realized by Dehn twists.

\thm{Fact}
{\label{fact2}
Let $\psi: F_2\rightarrow \bbbz^2$ be the canonical abelianizing
homomorphism. Then, if $\gamma_1$ and $\gamma_2$ are two primitive
elements in $F_2$ such that $\psi{\gamma_1}=\psi{\gamma_2},$ then
$\gamma_1$ and $\gamma_2$ are conjugate.
}

Fact \ref{fact2} follows easily from the work of Nielsen. For the
proof see \cite{osbzie}.

The following result, and its proof seem to go back to Poincar\'e.

\thm{Fact}
{\label{simgeod}
Let $S$ be a hyperbolic surface of finite volume, and let $\gamma$ be
a non-trivial simple closed curve, whose corresponding covering
transformation $\Gamma$ is hyperbolic. Then there is a unique
geodesic freely homotopic to $\gamma;$ this geodesic is simple.
}

\demo{Proof}
The existence of a unique geodesic $\tilde{\gamma}$ freely homotopic to
$\gamma$ follows by a completely standard straightening argument: since
$\Gamma$ is hyperbolic, it has two fixed points on the circle at
infinity of $\bbbh^2,$ as do all of its translates. The sought after
geodesic $\tilde{\gamma}$ is the unique geodesic in $\bbbh^2$ joining
those two fixed points. To show that $\tilde{\gamma}$ is simple is
also not hard. Indeed, the simplicity of $\gamma$ implies that for any
covering transformation $\beta,$ the fixed points of
$\beta\Gamma\beta^{-1}$ do not separate those of $\Gamma$ in the
cyclic order at infinity (otherwise the corresponding lifts of
$\gamma$ would have intersected.) However, this immediately implies
that the straightened curves do not intersect, by looking at it in the 
projective model of $\bbbh^2.$
\enddemo

\thm{Lemma}
{\label{homgeo}
No simple geodesic on the punctured torus separates.
}

\demo{Proof}
Suppose there was a separating simple geodesic $\gamma.$ Let it cut
the torus into two components $T_1$ and $T_2,$ and assume that $T_1$
is the component which does not contain the cusp. Double $T_1$ along
its boundary to obtain a compact oriented hyperbolic surface $S = 2
T_1.$ The area of $S$ is at least $4 \pi,$ and thus the area of $T_1$
is at least $2\pi.$ But the area of the whole punctured torus is
exactly $2\pi,$ and thus $\gamma$ could not have been separating.
\enddemo

\demo{Note}
The above argument can be easily adapted to show that there are no
simple closed geodesics on the thrice-punctured sphere.
\enddemo

The following lemma is well known:

\thm{Lemma}
{\label{bs0}
For any pair of non-separating simple closed curves $\gamma_1$ and
$\gamma_2$ on the punctured torus $T,$ there exists a homeomorphism of
$T$ taking $\gamma_1$ to $\gamma_2.$
}

\demo{Proof}
This is immediate from the classification of surfaces.
\enddemo

Fact \ref{fact2} and Lemmas \ref{homgeo} and \ref{bs0} combine to show
that is that there is a one-to-one correspondence between 
simple geodesics and primitive homology classes on
the punctured torus (a primitive homology class is a class $(m, n),$
such that $m$ and $n$ are relatively prime),
and it is geometrically obvious that this geodesic is really none other
than the familiar $(m, n)$ torus knot, or the $(m, n)$ geodesic on the
flat torus (without the puncture). The precise result is given
by the following construction of Osborne and Zieschang \cite{osbzie}:

Let $m$ and $n$ be a pair of relatively prime non-negative integers.
Define a function $f_{m, n}:\bbbz\rightarrow\{1,
2\}$ as follows: $f_{m, n}(k) = f_{m, n}(k'),$ if $k=k'$ modulo $m+n;$
for $k$ in $\{1, \dots, m\}$ let $f_{m, n}(k)=1,$ and for $k$ in $\{m+1,
\dots, m+n\}$ let $f_{m, n}(k) = 2.$

Now let $$W_{m, n}(x_1, x_2)=\prod_{i=0}^{m+n-1} x_{f_{m, n}(1+im)}.$$
If $m<0,$ let $W_{m, n}(x_1, x_2) = W_{-m, n}(x_1^{-1}, x_2),$ and if
$n<0,$ let $W_{m, n}(x_1, x_2) = W_{m, -n}(x_1, x_2^{-1}).$ 

The main theorem of \cite{osbzie} is the following:

\thm{Theorem}
{
If $x_1$ and $x_2$ are associated primitives in $F_2 = \langle s, t
\rangle,$ and if $m$ and $n$ are relatively prime, then $W_{m, n}(x_1,
x_2)$ is a primitive. furthermore, if $mq - np=1$ then $W_{m, n}$ and
$W_{p, q}$ are associated primitives. In particular, up to conjugation
the primitives of $F(s, t)$ are $\{W_{m, n}(s, t) | m, n \in \bbbz,
(m, n) = 1 \}.$
}

It can be seen that the words $W_{m, n}(s, t)$ are all cyclically
reduced, hence of minimal word length in their conjugacy class. Since
the word length of $W_{m, n}$ is equal to $m+n$ the following
combinatorial version of Theorem \ref{maintheorem} holds:

\thm{Theorem}
{
\label{combmain}
The number of conjugacy classes of primitive elements in $F(s, t)$
with reduced length not exceeding $L$ is asymptotic to $L^2.$
}

\demo{Proof} The number conjugacy classes in question equals
$$f(L)=|\{(m, n)| m+n\leq L, (m, n) = 1\}.$$ It follows from
elementary number theory that $f(L)$ is asymptotically equal to
$L^2/(2 \zeta(2)).$ 
\enddemo

\numbereddemo{Remark} For free groups of higher rank the analogue of
Theorem \ref{combmain} fails spectacularly, in that the number of
primitive elements of reduced length not exceeding $L$ grows
exponentially, as demonstrated by the following construction due to Casson
\cite{casson}. Let $F_n$ be the free group generated by $x_1, \dots,
x_n.$ Now let $y_i,$ $1\leq i \leq n$ be defined as follows: $y_i =
x_i$ when $i<n,$ and $y_n = x_n \prod x_{i_k}^{j_k},$ where $i_k<n$
for all $k, $ $i_k\neq i_{k+1}.$ Clearly, $\{y_i\}$ is a generating
set for $F_n, $ and also each $y_i$ is a cyclically reduced word. It
is also clear that the number of possibilities for $y_n$ of length not
exceeding $L$ grows exponentially in $L,$ for $n>2.$
\enddemo

\section{Geometry of group actions}
\label{cannon}

First, some definitions and a theorem, all
directly from \cite{canapp}. 

\demo{Definition}
A geometry is a metric space in which each bounded set has compact
closure. A group $G$ acts {\em geometrically} on a space $X$ if $X$ is a
geometry, and there is a homomorphism $\phi$ (usually suppressed) from
$G$ into the isometry group of $X$ such that the $G$-action is
properly discontinuous and cocompact ({\em properly discontinuous} means
that for each compact set $K$ in $X,$ the set $$\{g\in G | \emptyset
\neq  K \cap g K\}$$ is finite), and {\em cocompact} means that the
orbit space $X/G$ is compact).
\enddemo

\demo{Definition}
An intrinsic metric on a path space $X$ is one where the distance
between two points is the infimum of path lengths between those
points. 
\enddemo

\demo{Definition}
A relation $R: X\rightarrow Y$  between spaces $X$ and $Y$ is said to
be quasi-Lipschitz if $R$ is everywhere defined, and there exist
positive numbers $K$ and $L,$ such that for each $A\subset X,$ 
$$\diam R(A)\leq K \diam A + L.$$
\enddemo

\demo{Definition}
Relations $R: X\rightarrow Y$ and $S:Y\rightarrow X$ are {\em
quasi-inverses} if they are everywhere defined and there exists a
constant $M>0,$ such that $d(S \circ R, \id_X) < M$ an $d(R \circ S,
\id_Y) < M.$
\enddemo

\demo{Definition}
A relation $R: X\rightarrow Y$ is a {\em quasi-Lipschitz equivalence}
if there is a quasi-inverse $S:Y\rightarrow X,$ such that both $R$ and
$S$ are quasi-Lipschitz.
\enddemo

Finally, a theorem:

\thm{Theorem}
{\label{equiv}
If a group acts geometrically on two geometries $X$ and $Y$ with
intrinsic metrics, then $X$ and $Y$ are quasi-Lipschitz equivalent. 
}

Having dispensed with the generalities, an observation:

\thm{Observation}
{\label{ggeom}
If a finitely presented group $G$
is equipped with the word metric, then $G$ is a geometry,
on which $G$ acts geometrically by multiplication from the
right. 
}

Now, restrict to the case of current interest, where $G=F_2$ acts on
the hyperbolic plane $\bbbh^2$ with quotient the punctured torus. If
this action was geometric, in the sense of the above definitions, then
Theorem \ref{maintheorem} would follow from Theorem \ref{combmain},
Theorem \ref{equiv}, and the Observation \ref{ggeom}. However,
the action is not geometric, since the hyperbolic punctured torus is
not compact. It is not hard to construct the right geometry.

Consider the tessellation of
$\bbbh^2$ by fundamental domains for the given $F_2$ action, and
truncate each fundamental domain by a closed horocycle of length $2$
(as is possible by Theorem \ref{mcsimple}). The geometry $X$ will be
the subset of $\bbbh^2$ equal to the union of all these truncated
fundamental domains, and equipped with path metric. Clearly $X$ is a
geometry, on which $F_2$ acts geometrically, and so $X$ is
quasi-isometric (quasi-Lipschitz equivalent) to $F_2$ equipped with
its word metric.

While the general minimal paths on $X$ may be quite far from
hyperbolic geodesics, by Theorem \ref{mcsimple}, the 
minimal paths corresponding to lifts of simple closed geodesics are
unaffected by the truncation, and so their length are within a
constant factor of the length of the corresponding reduced words in
the word metric.

\section{Geometry of multicurves}
\label{multic}
Let $S$ be a hyperbolic surface of finite volume with at most one
cusp. We define a {\em multicurve} $m$ on $S$ to be a map from a (not
necessarily connected) 
$1$-manifold $M$ to $S.$ We define the {\em length} of $m$ to be the
sum of the lengths of the images of components of $M.$ We say that a
multicurve $m$ is {\em embedded} if the image of $m$ is the union of
simple closed curves $\gamma_1, \dots, \gamma_k$ on $S.$ Note that the
map $m$ may cover some of the components $\gamma_i$ multiple times, so
this does not coincide with the usual meaning of embedding. A
multicurve defines a singular chain, which, in turn, defines a
homology class in $H_1(S, \bbbz).$

The first observation is the following:

\thm{Theorem}{\label{minimal}
Let $h \in H_1(S, \bbbz)$ be a non-trivial homology class. There
exists a multicurve $m$ representing $h$ of minimal length, and $m$ is
embedded, with all components geodesic. 
}

\demo{Proof}
First, we show the existence. If $S$ is a compact surface, this follows by a
standard Arzela-Ascoli argument. If $S$ has cusps, let $l_h$ be the
infimum of the lengths of multicurves representing $h;$ clearly $l_h >
0.$ Let $m_1, \dots, m_l, \dots$ be multicurves whose lengths approach
$l_h.$ In order to apply the Arzela-Ascoli theorem, it is enough to
show that the diameters of the images of the $m_i$ are uniformly
bounded. Since the lengths of $m_i$ approach $l_h,$ it follows that
the diameters of all the components of $m_i,$ for $i$ sufficiently
great, are uniformly bounded by $2 l_h.$ On the other hand, we can
assume that no component of $m_i,$ for $i$ sufficiently large, can be
contained entirely in a horodisk of area $2$ surrounding each cusp,
since such a component is homologically trivial, and thus 
we can delete such a component with no change to homology, and replace
$m_i$ by the resulting multicurve $m'_i.$ Hence, we can assume that
all of the $m_i$ are contained in a compact subset of $S,$ and thus
the existence follows.

Now let $m$ be a multicurve of minimal length representing $h.$ It is
clear that each component of $m$ is geodesic. Suppose that $m$ is not
embedded, thus two components $c_1$ and $c_2$ of $m$ intersect. We can
assume that $c_1$ and $c_2$ are not multiply covered (if they are, we
can split off one circle off each), and since they are geodesic, the
intersection is transverse. Let $O$ be an intersection of $c_1$ and
$c_2,$  and let $A_1OB_1$ and $A_2OB_2$ be small directed segments of
$c_1$ and $c_2,$  respectively, surrounding $O.$ By cutting and
pasting, we can replace these by $s_1 = A_1 O B_2$ and $s_2=A_2 O
B_1,$ without changing the homology, and then by smoothing $s_1$ and
$s_2$ at $O,$ we obtain a shorter multicurve than $m$ representing
$h,$ thus arriving at a contradiction.
\enddemo

\thm{Corollary}{\label{mintor}
Let $T$ be a punctured torus equipped with a hyperbolic structure.
Then, the shortest multicurve representing a non-trivial homology
class $h$ is a simple closed geodesic if $h$ is a primitive homology
class (that is, not a multiple of another class), and a multiply
covered geodesic otherwise. In addition, the shortest multicurve
representing $h$ is unique.}

\demo{Proof} Since there is exactly one hyperbolic geodesic in any
homotopy class, it is enough to observe that any two non-homotopic
curves on a punctured torus intersect. 
Theorem \ref{minimal} then implies that the shortest multicurve
representing $h$ has one component, which is multiply covered if $h$
is not primitive (if $h=(m, n)$, with $d={\rm gcd}(m, n),$ then the
shortest multicurve representing $h$ is covered $d$ times). Since
there is at most one simple closed geodesic in a homology class (see,
eg section \ref{nielsen}), the uniqueness follows.
\enddemo

\section{A norm on homology of the punctured torus}

First, let us define a valuation $\ell$ on $H_1(T, \bbbz),$ where
$\ell(h)$ is defined to be the length of the shortest multicurve
representing $h.$ The valuation of the trivial homology class is defined to
be $0.$ Corollary \ref{mintor} implies that 
\begin{equation}
\label{linn}
\ell(n h)=n \ell (h),
\end{equation}
and that 
\begin{equation}
\label{triangle}
\ell(h+g)\leq \ell(h)+\ell(g),
\end{equation}
where the inequality is
strict if $h$ and $g$ are not both multiples of the same homology
class. The latter inequality holds, because the union of the shortest
multicurves corresponding to $h$ and $g$ is not embedded, and hence, 
the shortest multicurve corresponding to $h+g$ is shorter then the
union. It follows that $\ell$ can be extended to the rational homology
$H_1(T, \bbbq),$ by linearity (equation \ref{linn}), and further, to
$H_1(T, \bbbr)$ by continuity, which follows from equations \ref{linn}
and \ref{triangle}. Since $\ell(0)=0,$ and by equations \ref{linn} and
\ref{triangle}, $\ell$ is a pseudo-norm on the two-dimensional vector
space $H_1(T, \bbbr).$ By the results of sections 1, 2 and 3, $\ell$
is actually a norm, since the results of those sections imply that $0<
c_1<\ell(h)/\|h\|_1 < c_2 < \infty,$ where 
$\|h\|_1 = |m|+ |n|$ for $h=(m, n).$ It follows that the unit ball of
the norm $\ell$ is a compact, convex figure $B_\ell$ in the plane, and
the number of simple geodesics of length not exceeding $L$ on the
torus $T$ is equal to the number of primitive lattice points in $L
B_\ell.$ Thus, Theorem \ref{maintheorem} follows.

\section{Further investigations.}

The geometry of the unit ball of the norm $\ell$ is very interesting;
the authors will discuss it in a future paper. In particular, the
study of the unit ball can be used to show that the error term in
Theorem \ref{maintheorem} is essentially sharp. See \cite{mcriv2}.

The methods of this paper can be extended without much difficulty to
define a norm on homology of surfaces of higher genus, by extending
the length of shortest multicurve valuation, and this can be used to
count minimal multicurves. It seems non-trivial to use this to count
simple geodesics on surfaces of higher genus. Using the theory of
train tracks and the methods of this note it can be shown (McShane and
Rivin, in preparation) that the number of geodesics of length not 
exceeding $L$ on a closed surface of genus $g$ with $c$ cusps
has order of growth $L^{6g-6+2c}.$

\end{document}